\documentclass[12pt,a4paper]{article}
\usepackage{amsmath,amssymb,amsfonts}
\def\Bbb{\mathbb}

\title{\bf On the moduli of a Dedekind sum}

\author{Kurt Girstmair}
\date{}

\makeatletter
\let\@@maketitle=\maketitle
\def\maketitle{\def\thispagestyle##1{\relax}\@@maketitle}
\makeatother
\newtheorem{theorem}{Theorem}

\newtheorem{lemma}{Lemma}

\def\BE{\begin{equation}}
\def\EE{\end{equation}}
\def\BD{\begin{displaymath}}
\def\ED{\end{displaymath}}
\def\BA{\begin{array}}
\def\EA{\end{array}}
\def\BEA{\begin{eqnarray*}}
\def\EEA{\end{eqnarray*}}
\def\BI{\bibitem}

\def\N{\Bbb N}
\def\Z{\Bbb Z}

\def\R{\Bbb R}

\def\phi{\varphi}

\def\MB{\mbox}

\def\OV{\overline}

\def\DIV{\,|\,}

\def\MN{\medskip\noindent}
\def\STOP{\hfill$\Box$}

\def\DED{Dedekind }

\begin{document}
\maketitle

\begin{abstract}

\noindent
Let $s(a,b)$ denote the classical Dedekind sum and $S(a,b)=12s(a,b)$. Let
$k/q$, $q\in \N$, $k\in \Z$, $(k,q)=1$, be the value of $S(a,b)$. In a previous paper we showed that
there are pairs $(a_r,b_r)$, $r\in\N$, such that $S(a_r,b_r)=k/q$ for all $r\in \N$, the $b_r$'s growing in $r$ exponentially.
Here we exhibit such a sequence with $b_r$ a polynomial of degree $4$ in $r$.

\end{abstract}

\section*{1. Introduction}

Let $a$ be an integer, $b$ a natural number, and $(a,b)=1$. The classical \DED sum $s(a,b)$ is defined by
\BD
   s(a,b)=\sum_{k=1}^{b} ((k/b))((ak/b)).
\ED
Here
\BD
  ((x))=\begin{cases}
                 x-\lfloor x\rfloor-1/2 & \MB{ if } x\in\R\smallsetminus \Z; \\
                 0 & \MB{ if } x\in \Z
        \end{cases}
\ED
(see \cite[p. 1]{RaGr}).
It is often more convenient to work with
\BD
 S(a,b)=12s(a,b)
\ED instead. We call $S(a,b)$ a {\em normalized} \DED sum.

Let $q$ be a natural number, $k$ an integer, $(k,q)=1$. If $a\in\Z$, $b\in \N$, are such that $S(a,b)=k/q$, we call the number $b$ a {\em modulus} for $k/q$
(observe that $S(a,b)=S(a',b)$ if $a\equiv a'$ mod $b$).

The case $q=1$ is trivial. Indeed, $S(a,b)\in \Z$ if, and only if, $b$ divides $a^2+1$. In this case $S(a,b)=0$ and the moduli for $k/q=0$ can be considered as known (see \cite[p. 28]{RaGr}).  
Accordingly, we assume $q\ge 2$ in what follows.

In the previous paper \cite{Gi4} we showed that for every value $k/q$ of a normalized \DED sum there are infinitely many moduli. The sequence $b_r$, $r\in\N$,
of moduli exhibited there grew in $r$ exponentially, however. The proof was based on a  periodic continued fraction attached to $k/q$. The exponential growth was a
consequence of the exponential growth of the denominators of the convergents of this continued fraction.

In the present paper we use a completely different technique in order to show that there is a substantially denser sequence $b_r$, $r\in \N$, of moduli for $k/q$. Indeed, $b_r$ is a polynomial of degree $4$ in $r$.
Nevertheless, we think that this result is still far from the truth since the sequence of {\em all} moduli seems to be much denser. Our sequence $b_r$ is described in detail in Theorem \ref{t4} below.

\section*{2. The details}

Our work is based on the following lemma, whose proof can be found  in the paper \cite{Gi3}.

\begin{lemma} 
\label{l1}
Let $b, q$ be natural numbers, $q\ge 2$, $a$ an integer, $(a,b)=1$. Then $S(a,b)$ takes the form $k/q$ for some $k\in\Z$, $(k,q)=1$ if,
and only if, $b$ has the form
\BD
\label{2.2}
 b=\frac{q(a^2+1)}t,
\ED
where $t$ is a natural number, $(t,q)=1$.
$($Note that the condition $(t,q)=1$ implies $t\DIV a^2+1$.$)$
\end{lemma} 

For the time being we assume that $S(a,b)=k/q$, $k\in\Z, q\in\N$, $(k,q)=1$. By Lemma \ref{l1} this means that
$b=q(a^2+1)/t$ for some $t\in\N$ with $(t,q)=1$.
The following lemma expresses $S(a,b)$ in terms of two normalized \DED sums.
For a proof, see \cite{Gi3}.

\begin{lemma} 
\label{l2}
Let $q\ge 2$ and $t$ be natural numbers, $(t,q)=1$, $a$ an integer, $(a,q)=1$, $t\DIV a^2+1$, and $b=q(a^2+1)/t$.
Then
\BD
 S(a,b)=-S(aq,t)+S(at^*,q)+\frac{(q^2-1)a}{tq},
\ED
where $t^*$ is an integer satisfying $tt^*\equiv 1 \mod q$.
\end{lemma} 

Now let $q\ge 2$, $r$ and $n$ be natural numbers, $n$ a divisor of $r^2q^2+1$. The fact that $-1$ is a quadratic residue mod $n$ restricts the possible choices of $n$.
Put
\BE
\label{2.3}
  t=\frac{r^2q^2+1}n.
\EE
Then $t$ is a natural number with $(t,q)=1$. Let $s$ be an integer with $(s,q)=1$. We define
\BD
  a=st-rq\enspace \MB{and} \enspace b=\frac{q(a^2+1)}t.
\ED
Since $a\equiv -rq\mod t$,we have $a^2+1\equiv r^2q^2+1\mod t$. In view of (\ref{2.3}), $b$ is a natural number.
The following theorem is fundamental for this work.

\begin{theorem} 
\label{t2}
In the above setting,
\BE
\label{2.1}
S(a,b)=S(rq^2,n)+S(s,q)-\frac{(q^2-1)r}{n}+sq-\frac sq.
 \EE
\end{theorem} 

\noindent
{\em Proof.}
From Lemma \ref{l2} we obtain
\BD
\label{2.4}
 S(a,b)=-S(aq,t)+S(at^*,q)+\frac{(q^2-1)a}{tq}.
\ED
Since $aq\equiv -rq^2\mod t$, we have $S(aq,t)=-S(rq^2,t)$. The reciprocity law for \DED sums (see \cite[p. 4, Th. 1]{RaGr})
says
\BD
S(rq^2,t)=-S(t,rq^2)+\frac{rq^2}t+\frac t{rq^2}+\frac 1{rq^2t}-3.
\ED
Furthermore,
\BD
S(t, rq^2)=S(\frac{r^2q^2+1}n,rq^2)=S((r^2q^2+1)n^*,rq^2),
\ED
where $n^*$ is an integer such that $nn^*\equiv 1\mod rq^2$ (see \cite[p. 26, formula (33c)]{RaGr}; note that $(n,rq)=1$, by (\ref{2.3})). Hence we have
\BD
S(t, rq^2)=S(n^*,rq^2)=S(n,rq^2).
\ED
Then the reciprocity law gives
\BD
S(n,rq^2)=-S(rq^2,n)+\frac n{rq^2}+\frac{rq^2}n+\frac 1{nrq^2}-3.
\ED
The normalized \DED sum
$S(at^*,q)$ can be evaluated as follows:
\BD
 S(at^*,q)=S((st-rq)t^*,q)=S(s-rqt^*,q)=S(s,q).
\ED
When we use $a=st-rq$ and (\ref{2.3}) together wit the above identities, we obtain
(\ref{2.1}).
\STOP

\bigskip
\noindent
Next we write $r_2=nr_1+r$, where $r_1$ is an integer, $r_1\ge 0$.
Put
\BD
  t_2=\frac{r_2^2q^2+1}n.
\ED
Since $r_2\equiv r$ mod $n$, $t_2$ is a natural number, $(t_2,q)=1$.
Put
\BD
 a_2=st_2-r_2q \enspace \MB{and}\enspace b_2=q\frac{a_2^2+1}{t_2}.
\ED
Again, $b_2$ is a natural number and Theorem \ref{t2} can be applied to this new situation.  It yields
\BD
\label{2.6}
 S(a_2,b_2)=S(r_2q^2,n)+S(s,q)-\frac{(q^2-1)r_2}{n}+sq-\frac sq.
\ED
Since $r_2\equiv r\mod n$, $S(r_2q^2,n)=S(rq^2,n)$.
Moreover, $(q^2-1)r_2/n=(q^2-1)r_1+(q^2-1)r/n$.
In view of (\ref{2.1}) we obtain
\BD
  S(a_2,b_2)=S(a,b)-(q^2-1)r_1.
\ED
As a consequence of Lemma \ref{l2}, one knows how to remove the summand $-(q^2-1)r_1$
(see \cite{Gi3}).
Indeed, put
\BD
\label{2.11}
  a_3=a_2+r_1t_2q,\enspace b_3=\frac{q(a_3^2+1)}{t_2}.
\ED
Then $S(a_3,b_3)=S(a,b)$.
So the value of $S(a_3,b_3)$ is independent of $r_1\ge 0$.

By means of some computation one sees
that $b_3$ is a polynomial in $r_1$ of degree 4 with leading coefficient $nq^5$.

Finally, we have to adapt the above context to the general situation when $a\in\Z$, $b\in \N$, $(a,b)=1$, are given and $S(a,b)=k/q$, $k\in\Z$, $q\in \N$, $q\ge 2$, $(k,q)=1$. By Lemma \ref{l1},
$b$ has the form $b=q(a^2+1)/t$ for a natural number $t$ dividing $a^2+1$, $(t,q)=1$. Since $t$ and $q$ are co-prime, there are integers $s$ and $r$ such that
$a=st-rq$. Because $(a,q)=1$, we have $(s,q)=1$. The number $r$ can be assumed positive since $a=(s+uq)t-(r+ut)q$ for an arbitrary natural number $u$. Now put $n=(r^2q^2+1)/t$. Since $-rq\equiv a$ mod $t$, $n$ is a natural number
and we have exactly the above situation. We collect our hitherto found results in the following theorem.

\begin{theorem} 
\label{t4}
Let $a\in\Z$, $b\in \N$ be given, $(a,b)=1$. Suppose that $S(a,b)=k/q$, $k\in\Z$, $q\in \N$, $q\ge 2$, $(k,q)=1$. Define $t\in\N$, $(t,q)=1$, by
\BD
  b=q\frac{a^2+1}t.
\ED
Let $s\in\Z$ and $r\in\N$ be such that
\BD
 a=st-rq.
\ED
Further, define  $n\in\N$ by
\BD
 n=\frac{r^2q^2+1}t.
\ED
Let $r_1\in\Z$ be $\ge 0$,
\BD
t_2=\frac{(r_1n+r)^2q^2+1}n,\enspace a_3=st_2+(r_1t_2-r_1n-r)q\enspace \MB{ and }\enspace b_3=q\frac{a_3^2+1}{t_2}.
\ED
Then
\BD
\label{2.12}
 S(a_3,b_3)=S(a,b).
\ED
Here $b_3$ is a polynomial in $r_1$ of degree 4 with leading coefficient $nq^5$.
\end{theorem} 

\medskip
\noindent
{\em Example.} Let $a=2$, $b=7$, so $S(a,b)=6/7$. We obtain $q=7$ and $t=5$. Then $a=st-rq$ with $s=6$ and $r=4$. Therefore $n=157$. For $r_1=1$, we have $t_2=8090$, $a_3=104043$ and $b_3=9366455$.
For $r_1=2$, we obtain $t_2=31561$, $a_3=628994$ and $b_3=87748619$. Finally, we note $a_3=1897961$, $b_3=358087303$ for $r_1=3$. In all cases $S(a_3,b_3)=6/7$, of course.
The ratio of $b_3$ to the value of the leading monomial $r_1^4nq^5$ of $b_3$ is $\approx 3.55, 2.08, 1.68$ for $r_1=1,2,3$, respectively. In the case $r_1=10$ this ratio is $\approx 1.18$.
Of course, it tends to $1$ if $r_1$ tends to infinity.

\MN
{\em Remarks.} 1. In accordance with \cite{Gi4}, one may use the continued fraction $[0,\OV{3,2,1}]=(\sqrt{37}-4)/7$ and its convergents $a/b$ of order $2+6m$, $m\ge 0$, which give $S(a,b)=6/7$. The first values of $a/b$ are
\BD
2/7, 302/1015, 44090/148183, 6436838/21633703, 939734258/3158372455.
\ED
So the first few terms of this exponentially growing sequence of $b$'s grow slower than those of the preceding example.

2. The sequence of {\em all} possible moduli of a given value $k/q$ of a \DED sum seems to be much denser than the sequences supplied by our methods.
In the case $k/q=6/7$ we know the pairs $(a,b)=(2,7)$, $(22,35)$, $(57, 182)$, $(128,203)$, $(107,350)$, $(50,427)$, $(72,595)$ and
$(235,742)$ that give $S(a,b)=6/7$.



\vspace{0.5cm}
\noindent
Kurt Girstmair            \\
Institut f\"ur Mathematik \\
Universit\"at Innsbruck   \\
Technikerstr. 13/7        \\
A-6020 Innsbruck, Austria \\
Kurt.Girstmair@uibk.ac.at

\end{document}